\documentclass[11pt,twoside,reqno]{amsart}
\usepackage{amsmath,amssymb}
\usepackage[colorlinks=true,linkcolor=blue,citecolor=blue,urlcolor=blue]{hyperref}
\newtheorem{theorem}{Theorem}[section]
\newtheorem{lemma}[theorem]{Lemma}
\newtheorem{proposition}[theorem]{Proposition}
\newtheorem{corollary}[theorem]{Corollary}
\theoremstyle{remark}
\newtheorem{remark}[theorem]{Remark}
\numberwithin{equation}{section}
\title[Quasihyperbolic geodesics in the twice-punctured plane]{Quasihyperbolic Geodesics\\ in the Twice-Punctured Plane}
\author{Yali Du}
\address{(Yali Du) School of Mathematics and Statistics, Hubei University, Wuhan, 430062, P. R. China}
\email{Duyali\_123@126.com}
\author{Yike Ma}
\address{(Yike Ma) School of Mathematics and Statistics, Hubei University, Wuhan, 430062, P. R. China}
\email{squirrel12@163.com}
\author{Shengyou Wen}
\address{(Shengyou Wen) School of Mathematics and Statistics, Hubei University, Wuhan, 430062, P. R. China}
\email{sywen\_65@163.com}
\begin{document}
\begin{abstract}
We study quasihyperbolic distance and geodesic uniqueness in the plane punctured at $0$ and $2$. For endpoints on the perpendicular bisector of the punctures and on opposite sides of the real axis, the distance is the smaller of the length of the vertical segment and the least length of a feasible detour. For symmetric endpoints $1\pm it$ with $t>0$, the geodesic is unique precisely when $\operatorname{arsinh} t+\arctan t<\pi$. At equality there are exactly three geodesic images, and above the threshold there are exactly two. For arbitrary endpoints, we derive variational formulas for the distance in terms of admissible logarithmic spirals and their contacts with the bisector. The minimizing candidates in these formulas give all geodesic images and yield a length-comparison criterion for uniqueness.
\end{abstract}
\maketitle
\section{Introduction}
The quasihyperbolic metric, introduced by Gehring and Palka \cite{r3}, is the length metric with density reciprocal to the distance to the boundary. Even for a plane with one or two points removed, uniqueness of shortest curves depends on the positions of the endpoints and the available routes around the punctures. The map $f(w)=e^w$ provides a common starting point: in a once-punctured plane it converts quasihyperbolic length into Euclidean length, while in a Voronoi half-plane of a twice-punctured plane it does the same subject to an explicit admissibility constraint.

Existence of quasihyperbolic geodesics and their role in uniform-domain geometry were established by Gehring and Osgood \cite{r2}. Martin \cite{r7} proved local Lipschitz continuity of the first derivatives of their Euclidean arclength parametrizations; see also Väisälä \cite{r11}. Uniqueness has been studied under convexity assumptions and at small quasihyperbolic distances. Väisälä \cite{r10} treated convex domains in reflexive Banach spaces, and Martio and Väisälä \cite{r9} proved uniqueness in convex domains of uniformly convex Banach spaces. For planar domains, Väisälä \cite{r12} studied local uniqueness, prolongation and quasihyperbolic disks. Luiro \cite{r6} proved uniqueness in simply connected planar domains and the sharp sufficient condition $k_\Omega(x,y)<\pi$ in arbitrary planar domains.

Our first step recalls the classical formula of Martin and Osgood \cite{r8} in $\mathbb C\setminus\{0\}$. A shortest curve is the exponential image of a segment to a nearest endpoint lift. There is exactly one such curve unless the endpoints lie on opposite rays, when there are exactly two. This model also supplies the local pieces used throughout the twice-punctured problem.

For the twice-punctured plane, we normalize the punctures to $0$ and $2$ and write
\[\begin{gathered}G=\mathbb C\setminus\{0,2\},\qquad L=\{\operatorname{Re} z=1\},\\G_-=\{\operatorname{Re} z<1\}\setminus\{0\},\qquad G_+=\{\operatorname{Re} z>1\}\setminus\{2\}.\end{gathered}\]
The line $L$ is the perpendicular bisector of the segment joining $0$ and $2$. The sets $G_-$ and $G_+$ are the parts of $G$ lying to the left and right of $L$, respectively. Reflection across $L$ is a quasihyperbolic isometry. For endpoints on $L$, a unique geodesic must be fixed pointwise by reflection and hence must be the vertical segment. Folding a competing curve into the left half-plane and taking a logarithm reduces the problem to Euclidean length in the region $e^u\cos v\le1$. Its boundary graphs $u=-\log\cos v$ are strictly convex. Their tangents describe feasible spirals and show that a folded nonvertical minimizer has exactly one excursion from the bisector. This yields an attained detour minimum $J(a,b)$ and a necessary and sufficient uniqueness condition. For $1\pm it$, the condition becomes
\[\operatorname{arsinh}t+\arctan t<\pi.\]
The bisector calculation is also the central piece of the general endpoint problem. We prove that geodesics joining two points in the same open Voronoi half-plane stay in its closure: they are either one spiral, or two tangent spirals connected by a vertical segment. If one endpoint is on $L$, the first contact with $L$ separates off a single spiral tail. If the endpoints lie in opposite half-planes, the first and last contacts separate off two tails. The middle portion is exactly a geodesic between bisector points, already covered by the preceding calculation. Thus the three endpoint configurations have length formulas involving feasible single-spiral lengths and the same function $W(u,v)=k_G(1+iu,1+iv)$. For arbitrary endpoints these are variational formulas, rather than a single elementary coordinate inequality. Uniqueness means that their globally minimizing concatenations produce exactly one curve image.

Related work includes the point-to-line formula of Chen, Chen and Qian \cite{r1} and its applications to twice- and thrice-punctured planes. Herron \cite{r4} develops the general geometry of quasihyperbolic metric universal covers. The conformal-metric framework of Herron and Martin \cite{r5} gives another route to the regularity used here. Our arguments use tangent continuity at contacts and give the remaining length comparisons directly. Similarities preserve quasihyperbolic length, so normalization transfers all results to any pair of distinct punctures.

In Section 2, we present the once-punctured model. Section 3 develops the folded geometry and determines geodesics between bisector points. Section 4 uses the same excursion obstruction to classify same-half-plane geodesics and reduce arbitrary endpoints to the bisector problem. In each part the results are followed by their proofs and consequences.

\noindent\textit{Notation.} For a proper plane domain $\Omega$, set
\[\begin{gathered}\delta_\Omega(z)=\operatorname{dist}(z,\partial\Omega),\qquad \ell_{k,\Omega}(\gamma)=\int_\gamma\frac{|dz|}{\delta_\Omega(z)},\\ k_\Omega(p,q)=\inf_{\gamma:p\to q}\ell_{k,\Omega}(\gamma).\end{gathered}\]
A curve attaining this infimum is a quasihyperbolic geodesic. Uniqueness always means uniqueness of the curve image. Endpoints are distinct unless stated otherwise; for equal endpoints the constant curve is the unique minimizing image. A logarithmic spiral segment includes a radial segment or a circular arc as a special case.
\section{The once-punctured plane}

In this section, we determine when quasihyperbolic geodesics in the once-punctured plane are unique. The exponential covering map converts quasihyperbolic lengths into Euclidean lengths in logarithmic coordinates. Comparing the straight segments to the different lifts of the terminal point gives the classical distance formula and identifies exactly when two distinct geodesics occur.

\begin{theorem}\label{result:2.1}
Let $D=\mathbb C\setminus\{0\}$.
For distinct $x,y\in D$, the following conditions are equivalent:
\begin{enumerate}
\item
The quasihyperbolic geodesic from $x$ to $y$ is unique.
\item
The smaller angle between $x$ and $y$ is strictly less than $\pi$.
\item
There is no $\lambda>0$ such that $y=-\lambda x$.
\item
$0\notin[x,y]$.
\end{enumerate}
If they fail, there are exactly two geodesic images.
\end{theorem}

\begin{lemma}\label{result:2.2}
Let $f(w)=e^w$. For every absolutely continuous curve $\sigma:[a,b]\to\mathbb C$,
\begin{equation}\ell_{k,D}(f\circ\sigma)=\ell(\sigma)=\int_a^b|\sigma^\prime(t)|\,dt,\label{eq:2.1}\end{equation}
where $\ell(\sigma)$ denotes the Euclidean length of $\sigma$.
Conversely, an absolutely continuous curve in $D$ has a unique lift once its initial lift is prescribed. This lift is absolutely continuous and satisfies the same length identity.
\end{lemma}
\begin{proof}
Since $\delta_D(z)=|z|$ and $f^\prime=f$, the chain rule gives
\[\int_a^b\frac{|f(\sigma(t))\sigma^\prime(t)|}{|f(\sigma(t))|}\,dt=\int_a^b|\sigma^\prime(t)|\,dt.\]
The path-lifting property gives the unique lift. Locally it is the composition with a smooth inverse branch. A finite subdivision of the compact parameter interval suffices, with each inverse branch Lipschitz on a neighborhood of the corresponding compact image. Hence the lift is absolutely continuous, and the identity applies to it.
\end{proof}
Write $x=re^{i\alpha},\ y=se^{i\beta}$. Fix $A=\log r+i\alpha$; the terminal lifts are $B_n=\log s+i(\beta+2\pi n),\ n\in\mathbb Z$. The unique shortest lift to $B_n$ is $[A,B_n]$, with length
\begin{equation}L_n=\sqrt{\log^2(s/r)+(\beta-\alpha+2\pi n)^2}.\label{eq:2.2}\end{equation}
Thus the metric is the flat-cylinder metric $\mathbb R\times(\mathbb R/2\pi\mathbb Z)$.
\begin{proof}[Proof of Theorem \ref{result:2.1}]
Choose $\delta\in(-\pi,\pi]$ congruent to $\beta-\alpha$. If $|\delta|<\pi$, its absolute value is strictly smaller than that of every other representative, giving one nearest terminal lift. If $\delta=\pi$, exactly the changes $\pi,-\pi$ minimize \eqref{eq:2.2}; their projected curves are distinct. Finally an angle of $\pi$ is equivalent to $y=-\lambda x,\ \lambda>0$, and to $0\in[x,y]$.
\end{proof}

\begin{corollary}\label{result:2.3}
If $\vartheta\in[0,\pi]$ is the smaller angle between $x,y$, then
\[k_D(x,y)=\sqrt{\log^2\frac{|y|}{|x|}+\vartheta^2}.\]
For each minimizing angular change $\delta$, the corresponding curve is $\gamma(t)=\exp((1-t)A+tB)$, where $e^A=x$ and $B=A+\log(|y|/|x|)+i\delta$.
\end{corollary}
\begin{proof}
Minimize \eqref{eq:2.2} over the terminal lifts and apply Lemma \ref{result:2.2}. This is the classical formula of \cite{r8}.
\end{proof}
\begin{corollary}\label{result:2.4}
For distinct $x,y\in\mathbb C\setminus\{\xi\}$, the geodesic is unique if and only if $\xi\notin[x,y]$; otherwise there are exactly two.
\end{corollary}
\begin{proof}
Translate by $-\xi$ and apply Theorem \ref{result:2.1}.
\end{proof}
\section{The twice-punctured plane: bisector endpoints}

Throughout Sections 3--4 the domain and half-planes are those fixed in the introduction. We use $G_-$, $G_+$ for the open regions and $G_-\cup L$, $G_+\cup L$ for their closures in $G$. Parametrize $L$ by $z=1+iy$ with $y\in\mathbb R$. Then $dz=i\,dy$, so $|dz|=|dy|$, and $\delta_G(1+iy)=|1+iy|=|{-1}+iy|=\sqrt{1+y^2}$. Thus along $L$ we have
\begin{equation}\begin{gathered}\frac{|dz|}{\delta_G(z)}=\frac{|dy|}{\sqrt{1+y^2}},\qquad H(t)=\operatorname{arsinh}t,\\\ell_k([1+iu,1+iv])=|H(u)-H(v)|.\end{gathered}\label{eq:3.1}\end{equation}
For $t\in\mathbb R$ set $r_t=\sqrt{1+t^2}$, and for $c,d>0$ put
\[\Theta(c,d)=2\pi-\arctan c-\arctan d,\qquad M(c,d)=\frac{\log(r_c/r_d)}{\Theta(c,d)}.\]
For $a,b>0$ define
\begin{equation}\begin{split}J(a,b)=\min_{\substack{c\ge a,\ d\ge b\\-c\le M(c,d)\le d}}\big\{&H(c)-H(a)+H(d)-H(b)\\&+\sqrt{\log^2(r_c/r_d)+\Theta(c,d)^2}\big\}.\end{split}\label{eq:3.2}\end{equation}
Attainment is proved below.
\begin{theorem}\label{result:3.1}
Let $p=1+iy_1,\ q=1+iy_2$ be distinct points on $L$.
\begin{enumerate}
\item
If $y_1y_2\ge0$, the vertical segment is the unique geodesic.
\item
If $y_1y_2<0$, label the endpoints $p=1+ia,\ q=1-ib,\ a,b>0$. Then
\begin{equation}k_G(p,q)=\min\{H(a)+H(b),J(a,b)\}.\label{eq:3.3}\end{equation}
The geodesic is unique if and only if
\begin{equation}J(a,b)>H(a)+H(b).\label{eq:3.4}\end{equation}
When unique, it is vertical. Equality or the reverse strict inequality implies nonuniqueness.
\end{enumerate}
\end{theorem}
For later use, encode this result by
\begin{equation}W(u,v)=k_G(X_u,X_v),\qquad X_u=1+iu.\label{eq:3.5}\end{equation}
Thus $W(u,u)=0$, and
\begin{equation}W(u,v)=\begin{cases}|H(u)-H(v)|,&uv\ge0,\\\min\{H(a)+H(b),J(a,b)\},&uv<0,\end{cases}\label{eq:3.6}\end{equation}
where the positive and negative heights in the second line are $a,-b$. This is the middle-segment distance used in Section 4.

Reflection and folding across $L$ are
\[R(z)=2-\overline z,\qquad F(x+iy)=1-|x-1|+iy.\]
The fold preserves the density and Euclidean speed almost everywhere along absolutely continuous curves. On the set where the real coordinate equals 1, its derivative is zero almost everywhere, which also justifies the speed assertion there. Consequently
\begin{equation}\ell_k(F\circ\gamma)=\ell_k(\gamma).\label{eq:3.7}\end{equation}
Whenever the endpoints are fixed by the fold, a folded global minimizer is again globally minimizing. This applies both to endpoints on $L$ and to endpoints in $G_-\cup L$, as needed in Section 4.

On $G_-\cup L$, the density is $1/|z|$. Its full logarithmic preimage is
\[E=\{u+iv:e^u\cos v\le1\}.\]
Lengths in $E$ are Euclidean by Lemma \ref{result:2.2}. Define $g(v)=-\log\cos v$ for $|v|<\pi/2$. Then $g^\prime(v)=\tan v,\ g^{\prime\prime}(v)=\sec^2 v>0$, and within $|v-2\pi j|<\pi/2$ the constraint is
\begin{equation}u\le g(v-2\pi j).\label{eq:3.8}\end{equation}
When $\cos v\le0$, there is no restriction on $u$. In slope arguments we put $v$ horizontally and $u$ vertically; a slope means $du/dv$.

Global minimizers exist by \cite[Lemma~1]{r2}. Directly in $G$, the inequalities $|d\log|z||\le \frac{|dz|}{\delta_G(z)}$ and $|d\log|z-2||\le \frac{|dz|}{\delta_G(z)}$ keep a bounded-length minimizing sequence in a compact subset of $G$. Euclidean arclength bounds, compactness and lower semicontinuity yield a minimizer. Their Euclidean arclength parametrizations have locally Lipschitz first derivatives by \cite[Theorems~2.2 and~4.3]{r7}; see also \cite[Theorem~2.8]{r11}. In particular their tangents are continuous, including for folded curves when these are global minimizers. Alternatively, $\log(1/\delta_G)=\max\{-\log|z|,-\log|z-2|\}$ is locally Lipschitz and subharmonic, and the metric is complete by \cite[Corollary~2.2]{r3}; the regularity theorem of \cite{r5} applies with $\phi(t)=e^t$.

An excursion of a curve in $G_-\cup L$ is a nonconstant maximal portion with interior in $G_-$ and both endpoints on $L$. For a geodesic its logarithmic lift is a straight segment: local Euclidean straightness on overlapping intervals gives a fixed direction.
\begin{lemma}\label{result:3.2}
Every geodesic excursion has one positive-height endpoint and one negative-height endpoint. Up to reversal and a deck translation its lift joins
\[(\log r_c,\arctan c)\quad\text{to}\quad(\log r_d,2\pi-\arctan d),\qquad c,d>0.\]
\end{lemma}
\begin{proof}
A lift of $1+iy$ has coordinates $(g(\alpha),\alpha+2\pi j)$, where $\alpha=\arctan y$. A chord between distinct points on the same strictly convex boundary graph lies above it and is infeasible; equal endpoint lifts give only the constant segment. If the boundary indices differ by at least two, an intermediate angle $2\pi j$ lies between the endpoint angles. The linear height there is nonnegative, whereas feasibility requires $u\le0$. A positive value is impossible and a zero value gives an interior contact at $z=1$, also impossible for an excursion. The indices therefore differ by one.
Orient the lift from $(g(\alpha),\alpha)$ to $(g(\beta),\beta+2\pi)$, and let
\[m=\frac{g(\beta)-g(\alpha)}{2\pi+\beta-\alpha}.\]
Endpoint feasibility implies
\begin{equation}\tan\beta\le m\le\tan\alpha.\label{eq:3.9}\end{equation}
Thus $\beta\le\alpha$. If $0\le\beta\le\alpha$, monotonicity of $g$ gives $m\le0$, while \eqref{eq:3.9} gives $m\ge0$. Only $\alpha=\beta=0$ is possible, and that segment meets $L$ internally. The case $\beta\le\alpha\le0$ is excluded similarly. Hence $\beta<0<\alpha$; put $c=\tan\alpha,\ d=-\tan\beta$.
\end{proof}
\begin{remark}\label{result:3.3}
A one-turn excursion connects neighboring boundary components of $E$. Its actual angular change is $\Theta(c,d)\in(\pi,2\pi)$, not $2\pi$.
\end{remark}
\begin{lemma}\label{result:3.4}
A global geodesic contained in $G_-\cup L$ has no excursion whose two endpoints occur at interior parameter values of the entire curve.
\end{lemma}
\begin{proof}
Since the entire curve stays on one side of $L$, tangent continuity forces tangency to $L$ at both contacts. By Lemma \ref{result:3.2}, the excursion endpoints have heights $c,-d,\quad c,d>0$. Tangency of its straight logarithmic lift at these two endpoints would require its slope to be both $c$ and $-d$, which is impossible. The argument also applies to isolated contacts, without adjacent vertical intervals.
\end{proof}
This contact obstruction controls both the bisector problem below and the same-half-plane problem in Section 4.
\begin{proposition}\label{result:3.5}
If $y_1y_2\ge0$, the vertical segment between $X_{y_1},X_{y_2}$ is the unique geodesic.
\end{proposition}
\begin{proof}
For nonnegative heights, fold a minimizer across $L$ and then across the real axis into the upper half-plane. Both folds preserve length and fix the endpoints. Lemma \ref{result:3.2} excludes excursions in $\{\operatorname{Re} z\le1,\ \Im z\ge0\}$, so the result lies on $L$. The real-axis fold preserves the real coordinate, and $F^{-1}(L)=L$; therefore the original curve also lies on $L$. No backtracking is possible. Negative heights follow by reflection, including zero heights.
\end{proof}
\begin{proposition}\label{result:3.6}
Let $a,b>0$. Every nonvertical folded global geodesic from $X_a$ to $X_{-b}$ has exactly one excursion and has the form
\begin{equation}X_a\longrightarrow X_c\longrightarrow X_{-d}\longrightarrow X_{-b},\qquad c\ge a,\ d\ge b.\label{eq:3.10}\end{equation}
The first and last portions are vertical and, when nonconstant, are traversed upwards.
\end{proposition}
\begin{proof}
A geodesic through 1 must be vertical by Proposition \ref{result:3.5} applied to its two subcurves. A nonvertical one therefore avoids 1. An upper-to-lower excursion has increasing logarithmic angle by Lemma \ref{result:3.2}. At any adjacent vertical interval, tangent continuity then forces upward motion, since $\arctan y$ increases with $y$; reverse excursions force downward motion.

By Lemma \ref{result:3.4}, every excursion meets an endpoint of the whole curve, so there are at most two. If there were two, they would be the initial and final portions, separated by a vertical interval. Equal excursion directions would make that interval join opposite signs and pass through 1. Opposite directions would demand conflicting vertical orientations. A constant intervening interval is also impossible because the one-sided angular directions disagree. There is exactly one excursion. Before it the vertical portion cannot cross 1; the excursion therefore starts above and ends below the real axis. The upward orientations give $c\ge a,\ d\ge b$.
\end{proof}

\begin{lemma}\label{result:3.7}
The straight lift joining the points
\[(\log r_c,\arctan c)\qquad\text{and}\qquad(\log r_d,2\pi-\arctan d)\]
is feasible exactly when $-c\le M(c,d)\le d$. Its exponential image has interior strictly in $G_-$, and its length is
\[\sqrt{\log^2(r_c/r_d)+\Theta(c,d)^2}.\]
\end{lemma}
\begin{proof}
Its slope is $m=-M(c,d)$. Necessity follows from \eqref{eq:3.9}, namely $-d\le m\le c$. For sufficiency, near the initial boundary graph the line is below its supporting tangent because $m\le c$, and that tangent is below the strictly convex graph. Working backwards at the terminal graph gives the same conclusion from $m\ge-d$. Between angles $\pi/2$ and $3\pi/2$ there is no height restriction. Strict convexity gives strict interior separation from the graphs, and in the unrestricted interval $e^u\cos v\le0<1$. The length follows from Lemma \ref{result:2.2}.
\end{proof}
\begin{lemma}\label{result:3.8}
The minimum defining $J(a,b)$ is attained.
\end{lemma}
\begin{proof}
The feasible set is nonempty: take $c=d\ge\max\{a,b\}$. It is closed since $\Theta>\pi$ and all functions are continuous. The objective is continuous and bounded below by $H(c)+H(d)-H(a)-H(b)$, which tends to infinity as either variable tends to infinity while both have positive lower bounds. Nonempty bounded sublevel sets are compact, proving attainment.
\end{proof}
\begin{proof}[Proof of Theorem \ref{result:3.1}]
The same-sign statement is Proposition \ref{result:3.5}. For opposite signs the vertical length is $V=H(a)+H(b)$. Every feasible pair in \eqref{eq:3.2} gives an actual curve \eqref{eq:3.10} with the indicated objective length, so $k_G\le\min\{V,J\}$. Fold a global minimizer. Proposition \ref{result:3.6} says it is vertical or is one of those detours, giving the reverse inequality.

If $J>V$, all folded minimizers are vertical; $F^{-1}(L)=L$ forces the original ones to be vertical too. If $J=V$, an attained left detour and its reflection coexist with the vertical segment. If $J<V$, an attained detour and its reflection are distinct global minimizers. This proves both the formula and the uniqueness criterion.
\end{proof}

\begin{corollary}\label{result:3.9}
For $t>0$,
\[\begin{gathered}J(t,t)=2\pi-2\arctan t,\\ k_G(1+it,1-it)=\min\{2\operatorname{arsinh}t,2\pi-2\arctan t\}.\end{gathered}\]
Uniqueness is equivalent to $\operatorname{arsinh} t+\arctan t<\pi$. At equality there are exactly three geodesics, the vertical segment and two reflected circular detours; above the threshold there are exactly two, those detours.
\end{corollary}
\begin{proof}
Put $\alpha=\arctan t,\ r=\sqrt{1+t^2}$. Each feasible left detour has a lift joining $(\log r,\alpha)$ to $(\log r,2\pi-\alpha)$; its Euclidean length is at least $2\pi-2\alpha$, even if it includes boundary portions. The segment between these points is feasible by Lemma \ref{result:3.7} with $c=d=t$, and attains the bound. Equality forces precisely that segment. Its exponential image is the unique folded detour
\[z(v)=re^{iv},\qquad\alpha\le v\le2\pi-\alpha.\]
Its interior lies strictly on the left of $L$. A curve folding to it cannot switch sides internally without meeting $L$, so only this arc and its reflection are possible. Theorem \ref{result:3.1} now proves the counts and criterion.
\end{proof}
The function $\operatorname{arsinh} t+\arctan t$ starts at zero, tends to infinity, and has derivative $(1+t^2)^{-1/2}+(1+t^2)^{-1}>0$. Thus there is exactly one threshold,
\begin{equation}\operatorname{arsinh}t_0+\arctan t_0=\pi,\qquad t_0\approx3.1840793803.\label{eq:3.11}\end{equation}
\begin{remark}\label{result:3.10}
Luiro's sufficient bound $k_G(p,q)<\pi$ \cite[Theorem~1.2]{r6} is not necessary in this configuration. Since $\operatorname{arsinh}t_0=\pi-\arctan t_0>\pi/2$, we have $t_0>\sinh(\pi/2)$. For $\sinh(\pi/2)\le t<t_0$, the unique vertical geodesic has length $2\operatorname{arsinh}t\ge\pi$.
\end{remark}
The bisector distance $W$ in \eqref{eq:3.6} will now be used as the middle term in the general-endpoint length comparisons. The same folding map and Lemma \ref{result:3.4} also yield the required structure for endpoints away from $L$.
\section{The twice-punctured plane: arbitrary endpoints}
We now keep the domain, the open half-planes, and the functions $H$, $X_u$, $W$ fixed in Section 3. The role of $L$ changes: its first and last contacts with a geodesic separate the single-spiral tails from a middle segment whose distance is already known. For two endpoints in the same open half-plane, the contact obstruction of Lemma \ref{result:3.4} gives the stronger conclusion that the geodesic never enters the other half-plane.

\begin{proposition}\label{result:4.1}
Let $p,q\in G_-$. Every quasihyperbolic geodesic from $p$ to $q$ is contained in $G_-\cup L$. Its intersection with $L$, if nonempty, is connected. Consequently its image is either a single logarithmic spiral segment centered at 0, or has the form
\begin{equation}p\longrightarrow X_u\longrightarrow X_v\longrightarrow q,\label{eq:4.1}\end{equation}
where the first and last portions are logarithmic spiral segments in $G_-$, and the middle portion is the vertical segment in $L$. A one-point contact is allowed. If the vertical segment is nonconstant, both spiral portions are tangent to $L$ at their contacts. The reflected statement holds for $p,q\in G_+$, with center 2.
\end{proposition}
For use in the distance formulas we define a feasible single-spiral length. Let $c=0$ for the left half-plane and $c=2$ for the right. For endpoints in the corresponding closed punctured half-plane, choose arguments $\alpha_x=\arg(x-c),\ \alpha_y=\arg(y-c)$, and put
\begin{equation}\theta_{c,n}=\alpha_y-\alpha_x+2\pi n,\qquad n\in\mathbb Z,\label{eq:4.2}\end{equation}
\begin{equation}\begin{split}\Gamma_{c,n}(t)=c+\exp\big(& (1-t)\log|x-c|+t\log|y-c|\\ &+i(\alpha_x+t\theta_{c,n})\big),\end{split}\label{eq:4.3}\end{equation}
\begin{equation}S_{c,n}(x,y)=\sqrt{\log^2\frac{|y-c|}{|x-c|}+\theta_{c,n}^2}.\label{eq:4.4}\end{equation}
A branch is feasible when its entire image belongs to that closed half-plane. Define
\begin{equation}A_c(x,y)=\min_{n:\,\Gamma_{c,n}\ \mathrm{feasible}}S_{c,n}(x,y),\label{eq:4.5}\end{equation}
with value $+\infty$ if no branch is feasible and $A_c(x,x)=0$. The feasible condition is $\operatorname{Re}\Gamma_{0,n}(t)\le1$ or $\operatorname{Re}\Gamma_{2,n}(t)\ge1$, respectively, for every $t\in[0,1]$. Thus feasibility is part of the definition, not an additional assumption in the formulas below. A finite minimum in \eqref{eq:4.5} is attained, since $S_{c,n}\to\infty$ as $|n|\to\infty$. Changing the chosen arguments merely reindexes the integers.
\begin{theorem}\label{result:4.2}
For distinct endpoints the following formulas hold.
\begin{enumerate}
\item
If $p,q\in G_-$, then
\begin{equation}\begin{split} k_G(p,q)=\min\big\{A_0(p,q),\inf_{u,v\in\mathbb R}\big[&A_0(p,X_u)+|H(u)-H(v)|\\&+A_0(X_v,q)\big]\big\}.\end{split}\label{eq:4.6}\end{equation}
For $p,q\in G_+$, replace $A_0$ by $A_2$.
\item
If $p\in G_-$ and $q=X_b\in L$, then
\begin{equation}k_G(p,q)=\inf_{u\in\mathbb R}\big[A_0(p,X_u)+W(u,b)\big].\label{eq:4.7}\end{equation}
\item
If $p\in G_-$ and $q\in G_+$, then
\begin{equation}k_G(p,q)=\inf_{u,v\in\mathbb R}\big[A_0(p,X_u)+W(u,v)+A_2(X_v,q)\big].\label{eq:4.8}\end{equation}
\end{enumerate}
In \eqref{eq:4.7} and \eqref{eq:4.8} the infima are attained. In \eqref{eq:4.6}, whenever the second term supplies a global minimizer meeting $L$, its infimum is attained as well.
\end{theorem}

\begin{proof}[Proof of Proposition \ref{result:4.1}]
Apply Lemma \ref{result:3.4} to a folded global minimizer. Let $\gamma$ join $p,q\in G_-$, and fold it across $L$ using $F(x+iy)=1-|x-1|+iy$. The endpoints are fixed and length is preserved, so $F\circ\gamma$ is again a global minimizer. If $\gamma$ entered $G_+$, any component of its time in $G_+$ would have two interior endpoints on $L$. After folding it would be an excursion prohibited by Lemma \ref{result:3.4}. Thus $\gamma\subset G_-\cup L$.

If the contact set with $L$ were disconnected, there would be an excursion between two interior contact times, again impossible. The contact set is therefore an interval (possibly a singleton). The curve cannot backtrack, so its portion in $L$ is the vertical segment. Outside the contact interval its logarithmic lift is locally Euclidean straight, hence a single straight segment on each connected portion. Exponentiation gives the asserted spirals. Tangent continuity proves the contact condition for a nonconstant vertical portion. If the contact is a singleton, the two straight logarithmic pieces have the same direction and combine into one spiral segment. Reflection proves the right-half-plane assertion.
\end{proof}
\begin{proof}[Proof of Theorem \ref{result:4.2}]
On the closed left half-plane $\delta_G(z)=|z|$, and on the closed right half-plane $\delta_G(z)=|z-2|$. Lemma \ref{result:2.2}, after translation when necessary, therefore gives exactly \eqref{eq:4.4} for every feasible branch. A vertical segment from $X_u$ to $X_v$ has length $|H(u)-H(v)|$. Every finite expression in \eqref{eq:4.6} consequently is the length of an actual joining curve and is at least $k_G(p,q)$. Conversely, apply Proposition \ref{result:4.1} to a global minimizer. It is represented either by the direct term or by the three-piece term. Replacing any spiral portion by a shorter feasible branch cannot lower the global minimum, so its branch lengths equal the corresponding $A_0$ values. This proves \eqref{eq:4.6} and its attainment statement.

For \eqref{eq:4.7}, take a global minimizer and its first contact $X_u$ with $L$. Its initial portion lies in $G_-$, and a continuous logarithm makes it a Euclidean straight segment. Its length is therefore $A_0(p,X_u)$. The remaining subcurve is globally minimizing and has length $W(u,b)$. This gives a parameter attaining the distance. Conversely, each finite expression in \eqref{eq:4.7} can be realized by a feasible spiral followed by a geodesic between its two endpoints on $L$, and so cannot be smaller than $k_G(p,q)$.

For \eqref{eq:4.8}, every joining curve meets $L$. Choose the first and last contacts $X_u,X_v$ of a global minimizer. The initial and final portions lie in $G_-$ and $G_+$, respectively, and are single spirals. Its middle subcurve has length $W(u,v)$. The same replacement and concatenation arguments prove both inequalities and attainment. The case $u=v$, where the middle curve is constant, is included.
\end{proof}

To use the scalar distance formulas for uniqueness, retain all feasible branches attaining each finite $A_c$, and all geodesic images attaining $W$. At any parameter tuple attaining the global distance, concatenate these portions (or retain the direct branch in \eqref{eq:4.6}). Identify concatenations with the same curve image. Denote the resulting set of distinct images by $\mathcal M(p,q)$.
\begin{corollary}\label{result:4.3}
In each of the three cases, $\mathcal M(p,q)$ is exactly the set of quasihyperbolic geodesic images joining $p$ to $q$. In particular,
\begin{equation}\text{the geodesic is unique}\quad\Longleftrightarrow\quad\operatorname{card}\mathcal M(p,q)=1.\label{eq:4.9}\end{equation}
Here $\operatorname{card}\mathcal M(p,q)$ denotes the number of distinct quasihyperbolic geodesic images joining $p$ to $q$. Different parametrizations of the same image are counted only once.
Equivalently, among the feasible candidate images described by the corresponding formula, one image has strictly smaller length than every different image.
\end{corollary}
\begin{proof}
Each attaining concatenation has length $k_G(p,q)$, hence is a global geodesic. Conversely every global geodesic has the decomposition used in the proof of Theorem \ref{result:4.2}, with each constituent attaining its indicated minimum. Thus no geodesic image is omitted.
\end{proof}
In parts 2 and 3 of Theorem \ref{result:4.2}, every middle segment used in an attaining concatenation must itself be unique if the entire geodesic is unique. Theorem \ref{result:3.1} makes this condition explicit: same-side heights give a unique vertical segment; for heights $a,-b$, with $a,b>0$, uniqueness requires
\begin{equation}J(a,b)>H(a)+H(b).\label{eq:4.10}\end{equation}
This is necessary for the selected middle segment but does not by itself exclude a different minimizing contact configuration or a different minimizing tail branch.

For part 1 of Theorem \ref{result:4.2}, the nonconstant vertical candidates can be located through tangency equations. If $p=re^{i\alpha}$ and the initial spiral meets $X_u$ tangentially, its logarithmic straight segment satisfies
\begin{equation}\log\frac{\sqrt{1+u^2}}{r}=u(\arctan u-\alpha+2\pi n).\label{eq:4.11}\end{equation}
The analogous equation holds for $q$ and $X_v$. Feasibility and agreement of the oriented tangents with the vertical portion must also be imposed. Equation \eqref{eq:4.11} alone does not guarantee a minimizing curve.
\begin{remark}\label{result:4.4}
These results reduce the arbitrary-endpoint problem to explicit spiral lengths, one or two real contact heights, integer argument branches, and the function $W$ already determined in Section 3. They do not assert a single elementary coordinate inequality or a finite number of argument branches independent of the endpoints. Minimizing parameter tuples need not be unique representations of a curve: uniqueness is always counted by images. In particular, $A_c$ is a feasible single-spiral length, not the unrestricted intrinsic distance in a punctured half-plane.
\end{remark}

Theorem \ref{result:4.2} shows the precise difference between endpoint configurations. Same-half-plane endpoints require direct spirals and spiral--vertical--spiral paths only. A bisector endpoint requires one contact height and the function $W$. Opposite-half-plane endpoints require two heights and the same $W$; equal heights are allowed. These reductions classify possible minimizing images but do not imply that a minimizing pair of contact heights alone guarantees uniqueness.
\begin{corollary}\label{result:4.5}
Let $\xi_1,\xi_2\in\mathbb C$ be distinct and put $\Omega=\mathbb C\setminus\{\xi_1,\xi_2\}$. The similarity
\[T(z)=\frac{2(z-\xi_1)}{\xi_2-\xi_1}\]
induces a length-preserving bijection between joining curves in $\Omega$ and $G$. For arbitrary $p,q\in\Omega$,
\[k_\Omega(p,q)=k_G(T(p),T(q)),\]
and geodesic images correspond bijectively. Thus Theorem \ref{result:3.1} applies to endpoints on the perpendicular bisector of $\xi_1,\xi_2$, while Theorem \ref{result:4.2} and the length-comparison criterion apply to all other endpoint configurations.
\end{corollary}
\begin{proof}
The map sends the punctures to 0, 2, their bisector to $L$, and multiplies both Euclidean length elements and boundary distances by $2/|\xi_2-\xi_1|$. Their quotient is unchanged. Apply the preceding results after this normalization.
\end{proof}
\begin{remark}\label{result:4.6}
For $p=1/2,\ q=3/2$, the unique geodesic is the real segment and has length $2\log 2$. Indeed every joining curve meets $X_u=1+iu$. Using $\delta_G(z)\le|z|$ on its initial portion and $\delta_G(z)\le|z-2|$ on its final portion gives total length at least
\[2\log\big(2\sqrt{1+u^2}\big)\ge2\log 2.\]
The real segment attains this bound. Equality forces $u=0$ and both portions to be radial. In \eqref{eq:4.8} this is the attaining choice $u=v=0,\ W(0,0)=0$, with two radial tails of length $\log 2$.
\end{remark}

\end{document}